# A STOCHASTIC LOG-LAPLACE EQUATION[1]

### By Jie Xiong

#### *University of Tennessee*

We study a nonlinear stochastic partial differential equation whose solution is the conditional log-Laplace functional of a superprocess in a random environment. We establish its existence and uniqueness by smoothing out the nonlinear term and making use of the particle system representation developed by Kurtz and Xiong [*Stochastic Process. Appl.* **83** (1999) 103–126]. We also derive the Wong–Zakai type approximation for this equation. As an application, we give a direct proof of the moment formulas of Skoulakis and Adler [*Ann. Appl. Probab.* **11** (2001) 488–543].

## 1. Introduction and main results.

1.1. *Introduction.* We study the behavior of a branching interacting particle system in a random environment. For simplicity of notation, we assume that the particles move in the one-dimensional space $\mathbb{R}$. The branching is critical binary; that is, at independent exponential times, each particle will die or split into two with equal probabilities. Between branchings, the motion of the $i$th particle is governed by an individual Brownian motion $B_i(t)$ and a common Brownian motion $W(t)$ which applies to all particles in the system:

$$(1.1) \qquad d\eta_t^i = b(\eta_t^i)\,dt + c(\eta_t^i)\,dW(t) + e(\eta_t^i)\,dB_i(t), \qquad i = 1, 2, \ldots,$$

where $b$, $c$, $e$ are real functions on $\mathbb{R}$ ($c, e \geq 0$), $W$, $B_1$, $B_2$, ... are independent (standard) Brownian motions and $\eta_t^i$ is the position of the $i$th particle at time $t$. Let $\mathcal{M}_F(\mathbb{R})$ denote the set of all finite Borel measures on $\mathbb{R}$. It is established by Skoulakis and Adler [18] that the high-density limit $X_t$ of this system is the unique $\mathcal{M}_F(\mathbb{R})$-valued solution to the following *martingale*

Received May 2002; revised June 2003.

[1]Supported in part by the NSA.

*problem* (MP): $X_t$ is a continuous process with initial $X_0 = \mu \in \mathcal{M}_F(\mathbb{R})$ such that, for any $\phi \in C_b^2(\mathbb{R})$,

$$M_t(\phi) \equiv \langle X_t, \phi \rangle - \langle \mu, \phi \rangle - \int_0^t \langle X_s, b\phi' + a\phi'' \rangle \, ds$$

is a continuous martingale with quadratic variation process

$$\langle M(\phi) \rangle_t = \int_0^t (\langle X_s, \phi^2 \rangle + |\langle X_s, c\phi' \rangle|^2) \, ds,$$

where $a(x) = \frac{1}{2}(e(x)^2 + c(x)^2)$. Moment formulas are derived in [18]. A related model is studied by Wang [19] and Dawson, Li and Wang [4].

The log-Laplace equation has been used by many authors in deriving various properties for superprocesses (cf. [2, 5]). It is natural, as indicated in [18], to derive properties of $X_t$ by making use of the corresponding backward stochastic log-Laplace equation (LLE):

$$(1.2) \quad \begin{aligned} y_{s,t}(x) = f(x) &+ \int_s^t (b(x)\partial_x y_{r,t}(x) + a(x)\partial_x^2 y_{r,t}(x) - y_{r,t}(x)^2) \, dr \\ &+ \int_s^t c(x)\partial_x y_{r,t}(x) \, \hat{d}W_r, \end{aligned}$$

where $f$ is the test function for the Laplace transform [cf. (1.8)], $\partial_x$, $\partial_x^2$ are the first and second partial derivatives with respect to $x$ and the last integral is the backward Itô integral. Since a solution to (1.2) is not established in [18], the moment formulas for $X_t$ are derived based on other techniques. The establishment of a unique solution to (1.2) is posed by [18] as an interesting challenge.

In this paper, we study the LLE (1.2). The main result is Theorem 1.2 in which we prove that the log-Laplace transform of $X_t$ is indeed given by the solution to (1.2). For simplicity of notation, we consider the forward version of the LLE:

$$(1.3) \quad \begin{aligned} y_t(x) = f(x) &+ \int_0^t (b(x)\partial_x y_r(x) + a(x)\partial_x^2 y_r(x) - y_r(x)^2) \, dr \\ &+ \int_0^t c(x)\partial_x y_r(x) \, dW_r. \end{aligned}$$

The stochastic partial differential equation (SPDE) is an important field of current research. We refer the reader to [1, 11] and [17] for an introduction to this topic. Many authors studied linear SPDEs. Here we only mention two recent papers: [9] and [14]. Fine properties of the solutions are established. Nonlinear SPDEs have also been studied. Here we mention a sequence of papers by Kotelenez [12, 13] which are the closest to the present setting. In this case, the derivative of the solution is not involved in the noise term. To the best of our knowledge, the LLE (1.3) does not fit into the setups of existing theory of SPDE.

1.2. *Main results.* First we study the existence and uniqueness for the solution to (1.3). We also establish its particle system representation in the spirit of Kurtz and Xiong [15].

To begin with, we introduce some notation needed in this paper. Let $H_0 = L^2(\mathbb{R})$ be the set of all square integrable functions on $\mathbb{R}$, and let $H_0^+$ consist of all the nonnegative functions in $H_0$. Let $H_m = \{\phi \in H_0 : \phi', \ldots, \phi^{(m)} \in H_0\}$. Define the Sobolev norm on $H_m$ by

$$\|\phi\|_m^2 = \sum_{j=1}^{m} \int |\phi^{(j)}(x)|^2 \, dx.$$

Use $\langle \cdot, \cdot \rangle$ to denote the inner product in $H_0$ or the integral of a function with respect to a measure.

DEFINITION 1.1. An $H_0^+$-valued (measurable) process $y_t$ is a solution to (1.3) if, for any $\phi \in C_0^\infty(\mathbb{R})$,

$$\langle y_t, \phi \rangle = \langle f, \phi \rangle + \int_0^t \langle y_r, -(b\phi)' + (a\phi)'' - y_r\phi \rangle \, dr$$

$$+ \int_0^t \langle y_r, -(c\phi)' \rangle \, dW_r, \qquad t \geq 0.$$

Throughout this paper, we assume the following

BOUNDEDNESS CONDITION (BC). $f \geq 0, b, c, e$ are bounded functions with bounded first and second derivatives. Denote a bound by $K$. Further, $e$ is bounded away from 0, $c$ has third continuous and bounded derivative and $f$ is of compact support.

THEOREM 1.2. *Suppose that Condition* (BC) *holds. Then:*

  (i) *The LLE* (1.3) *has a unique solution* $y_t(x)$.
  (ii) $y_t$ *is the unique solution of the following infinite particle system:* $i = 1, 2, \ldots,$

$$d\xi_t^i = e(\xi_t^i) \, dB_i(t) + (2a' - b - cc')(\xi_t^i) \, dt - c(\xi_t^i) \, dW_t, \tag{1.4}$$

$$dm_t^i = m_t^i((a'' - b' - y_t)(\xi_t^i) \, dt - c'(\xi_t^i) \, dW_t), \tag{1.5}$$

$$Y_t = \lim_{n\to\infty} \frac{1}{n} \sum_{i=1}^{n} m_t^i \delta_{\xi_t^i} \qquad a.s. \tag{1.6}$$

*where, for any* $t \geq 0$, $Y_t$ *is absolutely continuous with respect to Lebesgue measure and* $y_t$ *is the Radon–Nikodym derivative.*

Next, we consider the Wong–Zakai type approximation to LLE (1.3):

$$y_t^\varepsilon(x) = f(x) + \int_0^t (\bar{b}(x)\partial_x y_r^\varepsilon(x) + \bar{a}(x)\partial_x^2 y_r^\varepsilon(x) - y_r^\varepsilon(x)^2)\, dr$$
(1.7)
$$+ \int_0^t c(x)\partial_x y_r^\varepsilon(x)\dot{W}_r^\varepsilon\, dr,$$

where $\bar{b}(x) = b(x) - \frac{1}{2}c(x)c'(x)$, $\bar{a}(x) = \frac{1}{2}e(x)^2$ and, for $k\varepsilon \le r < (k+1)\varepsilon$, $\dot{W}_r^\varepsilon = \varepsilon^{-1}(W_{(k+1)\varepsilon} - W_{k\varepsilon})$.

THEOREM 1.3. *Suppose that Condition* (BC) *holds. Then for any $t \ge 0$,*

$$\mathbb{E}\int |y_t^\varepsilon(x) - y_t(x)|^2\, dx \to 0$$

*as $\varepsilon \to 0$.*

Now we consider the Wong–Zakai approximation to the measure-valued process $X$. Let $\mathbb{P}^W$ be the conditional probability measure given $W$. Let $X^\varepsilon$ be the solution to the following *conditional martingale problem* (CMP): $X^\varepsilon$ is a continuous $\mathcal{M}_F(\mathbb{R})$-valued process such that, for any $\phi \in C_b^2(\mathbb{R})$,

$$M_t^\varepsilon(\phi) \equiv \langle X_t^\varepsilon, \phi\rangle - \langle X_0^\varepsilon, \phi\rangle - \int_0^t \langle X_s^\varepsilon, (\bar{b} + c\dot{W}_s^\varepsilon)\phi' + \bar{a}\phi''\rangle\, ds$$

is a continuous $\mathbb{P}^W$-martingale with quadratic variation process

$$\langle M^\varepsilon(\phi)\rangle_t = \int_0^t \langle X_s^\varepsilon, \phi^2\rangle\, ds.$$

Let $\bar{\mathbb{R}} \equiv \mathbb{R} \cup \{\partial\}$ be the one-point compactification of $\mathbb{R}$. Denote by $\mathcal{M}_F(\bar{\mathbb{R}})$ the space of all finite measures on $\bar{\mathbb{R}}$ with the weak convergence topology. Note that $\mathcal{M}_F(\mathbb{R})$ can be regarded as a subset of $\mathcal{M}_F(\bar{\mathbb{R}})$ by extending the measure at $\partial$ as 0.

THEOREM 1.4. *As $\varepsilon \to 0$, if $X_0^\varepsilon \to \mu$ in $\mathcal{M}_F(\mathbb{R})$, then $X^\varepsilon \to X$ in $C([0,\infty), \mathcal{M}_F(\bar{\mathbb{R}}))$ in conditional law $\mathbb{P}^W$ for almost all $W$. As a consequence, we have*

(1.8)            $$\mathbb{E}^W \exp(-\langle X_t, f\rangle) = \exp(-\langle \mu, y_{0,t}\rangle)      a.s.$$

Finally, we derive the moment formulas of $X_t$. Note that these formulas have been obtained in [18] by a different method. Let $p(t, x, y)$ and $q(t, (x_1, x_2), (y_1, y_2))$ be the transition density functions of the Markov processes with generators

$$\mathcal{L}_1\phi(x) = b(x)\phi'(x) + a(x)\phi''(x)$$

and

$$\mathcal{L}_2 F(x_1, x_2) = b(x_1)\partial_{x_1} F + b(x_2)\partial_{x_2} F$$
$$+ a(x_1)\partial_{x_1}^2 F + a(x_2)\partial_{x_2}^2 F + c(x_1)c(x_2)\partial_{x_1}\partial_{x_2} F,$$

respectively.

THEOREM 1.5. *Suppose that Condition* (BC) *holds. For any bounded continuous function* $f$, *we have*

$$(1.9) \qquad \mathbb{E}(\langle X_t, f \rangle) = \iint f(y)p(t, x, y)\,dy\mu(dx)$$

*and*

$$\mathbb{E}(\langle X_t, f \rangle^2)$$
$$= \int_{R^4} f(y_1)f(y_2)q(t, (x_1, x_2), (y_1, y_2))\,dy_1\,dy_2\mu(dx_1)\mu(dx_2)$$
$$(1.10) \qquad + 2\iint\int_0^t \int p(t - s, x, y)$$
$$\times \iint f(z_1)f(z_2)q(s, (y, y), (z_1, z_2))\,dz_1\,dz_2\,dy\,ds\mu(dx).$$

We shall use $K$ with a subscript to denote a constant. If it will be quoted, the subscript will be the equation where it is defined. Otherwise, we shall use $K_1, K_2, \ldots$ in the proof of a proposition and the sequence starts over again in the proof of a new proposition. For example, $K_1$ may appear in the proofs of two different propositions to represent different constants.

Note that the Wong–Zakai approximation is not really needed to obtain the results in Theorems 1.4 and 1.5. An easier approach in deriving (1.8) is available. We refer the reader to [16] for the treatment of a related model which adds immigration structure to a branching interacting system studied in [4] and [19]. In this paper, we use the Wong–Zakai approximation because this is part of the conjecture in [18] and the main purpose of the current paper is to solve that conjecture. Furthermore, the Wong–Zakai approximation is of interest on its own.

**2. Stochastic log-Laplace equation.** In this section, we prove Theorem 1.2.

2.1. *Approximation.* To establish the existence of a nonnegative solution to (1.3), we smooth and truncate its nonlinear term and consider

$$^{\varepsilon}y_t(x) = f(x)$$

(2.1)        $+ \int_0^t (b(x)\partial_x {}^\varepsilon y_r(x) + a(x)\partial_x^2 {}^\varepsilon y_r(x) - (T_\varepsilon {}^\varepsilon y_r^\varepsilon(x)){}^\varepsilon y_r(x))\, dr$

              $+ \int_0^t c(x)\partial_x {}^\varepsilon y_r(x)\, dW_r,$

where $T_\varepsilon h(x) \equiv \int p_\varepsilon(x-z)h(z)\, dz$, $p_\varepsilon(x) = (2\pi\varepsilon)^{-1/2}\exp(-\frac{1}{2\varepsilon}x^2)$, ${}^\varepsilon y_r^\varepsilon(x) = \hat{{}^\varepsilon y_r}{}^\varepsilon y_r(x)$ and

$$\hat{{}^\varepsilon y_r} = \frac{\int {}^\varepsilon y_r(u)\, du \wedge \varepsilon^{-1}}{\int {}^\varepsilon y_r(u)\, du}$$

with the convention that $\frac{0}{0} = 0$.

LEMMA 2.1.  *Equation* (2.1) *has a unique solution.*

PROOF.  Consider the following infinite particle system: $i = 1, 2, \ldots,$

              $d\xi_t^i = e(\xi_t^i)\, dB_i(t) + (2a' - b - cc')(\xi_t^i)\, dt - c(\xi_t^i)\, dW_t,$

(2.2)        $dm_t^{\varepsilon,i} = m_t^{\varepsilon,i}((a'' - b' - T_\varepsilon {}^\varepsilon Y_t^\varepsilon)(\xi_t^i)\, dt - c'(\xi_t^i)\, dW_t),$

              ${}^\varepsilon Y_t = \lim_{n\to\infty} \frac{1}{n}\sum_{i=1}^n m_t^{\varepsilon,i}\delta_{\xi_t^i} \qquad \text{a.s.},$

where $\forall \nu \in M_+(\mathbb{R})\nu^\varepsilon \in M_+(\mathbb{R})$ is defined by $\nu^\varepsilon = \frac{\nu(\mathbb{R}) \wedge \varepsilon^{-1}}{\nu(\mathbb{R})}\nu$.

Now we show that the conditions of [15] are satisfied by the coefficients of the system (2.2). We only check those for

$$d_\varepsilon(x, \nu) \equiv -(T_\varepsilon \nu^\varepsilon)(x).$$

The verification for other coefficients is trivial.

Note that $p_\varepsilon(x) \leq (\sqrt{2\pi\varepsilon})^{-1}$ and

$$|\partial_x p_\varepsilon(x)| \leq \frac{1}{\sqrt{2\pi\varepsilon}}\sup_x e^{-x^2/2\varepsilon}\frac{|x|}{\sqrt{\varepsilon}} = \frac{1}{\sqrt{2\pi e}\varepsilon}.$$

Then

$$|d_\varepsilon(x, \nu)| = \left| \int p_\varepsilon(x-y)\nu^\varepsilon(dy) \right| \leq (\sqrt{2\pi\varepsilon})^{-1}.$$

Let

$$\mathbb{B}_1 = \{g \in C(\mathbb{R}) : |g(x)| \leq 1, \; |g(x) - g(y)| \leq |x - y| \,\forall x, y \in \mathbb{R}\}$$

and

$$\rho(\nu_1, \nu_2) = \sup_{g \in \mathbb{B}_1} |\langle \nu_1 - \nu_2, g\rangle|.$$

For $g \in \mathbb{B}_1$, we have

$$
\begin{aligned}
|\langle \nu_1^\varepsilon - \nu_2^\varepsilon, g \rangle| &\leq \frac{\nu_1(\mathbb{R}) \wedge \varepsilon^{-1}}{\nu_1(\mathbb{R})} |\langle \nu_1 - \nu_2, g \rangle| \\
&\quad + |\langle \nu_2, g \rangle| \left| \frac{\nu_1(\mathbb{R}) \wedge \varepsilon^{-1}}{\nu_1(\mathbb{R})} - \frac{\nu_2(\mathbb{R}) \wedge \varepsilon^{-1}}{\nu_2(\mathbb{R})} \right| \\
&\leq \rho(\nu_1, \nu_2) + |\langle \nu_1 - \nu_2, 1 \rangle| + |\nu_1(\mathbb{R}) \wedge \varepsilon^{-1} - \nu_2(\mathbb{R}) \wedge \varepsilon^{-1}| \\
&\leq 3\rho(\nu_1, \nu_2).
\end{aligned}
$$

Then

$$
\begin{aligned}
|d_\varepsilon(x_1, &\nu_1) - d_\varepsilon(x_2, \nu_2)| \\
&\leq \left| \int (p_\varepsilon(x_1 - y) - p_\varepsilon(x_2 - y)) \nu_1^\varepsilon(dy) \right| \\
&\quad + \left| \int p_\varepsilon(x_2 - y) \nu_1^\varepsilon(dy) - \int p_\varepsilon(x_2 - y) \nu_2^\varepsilon(dy) \right| \\
&\leq (\sqrt{2\pi e}\varepsilon^2)^{-1} |x_1 - x_2| + (\sqrt{2\pi}\varepsilon)^{-1} (e\varepsilon \wedge 1)^{-1/2} \rho(\nu_1^\varepsilon, \nu_2^\varepsilon) \\
&\leq K_1 \sqrt{|x_1 - x_2|^2 + \rho(\nu_1, \nu_2)^2}.
\end{aligned}
$$

By [15], $^\varepsilon Y_t$ is the unique solution to

$$
\langle {}^\varepsilon Y_t, \phi \rangle = \langle f, \phi \rangle + \int_0^t \langle {}^\varepsilon Y_r, (a\phi)'' - (b\phi)' - (T_\varepsilon {}^\varepsilon Y_r^\varepsilon)\phi \rangle \, dr - \int_0^t \langle {}^\varepsilon Y_r, (c\phi)' \rangle \, dW_r.
$$

Further, $^\varepsilon Y_t$ has density $^\varepsilon y_t$ which belongs to $H_0$. $\quad\square$

2.2. *Boundedness.* In this section, we establish a comparison result for SPDEs of the form (2.1). As a consequence, we obtain the boundedness of $^\varepsilon y_t$.

LEMMA 2.2. *For all $r, x$, we have*

$$
{}^\varepsilon y_r(x) \leq \|f\|_\infty \qquad a.s.,
$$

*where $\|f\|_\infty$ is the supremum of $f$.*

PROOF. Let $\tilde{m}_t^i$ be given by

$$
d\tilde{m}_t^i = \tilde{m}_t^i ((a'' - b')(\xi_t^i) \, dt - c'(\xi_t^i) \, dW_t)
$$

and let

$$
\tilde{Y}_t = \lim_{n \to \infty} \frac{1}{n} \sum_{i=1}^n \tilde{m}_t^i \delta_{\xi_t^i} \qquad \text{a.s.}
$$

Then $m_t^{\varepsilon,i} \leq \tilde{m}_t^i$ and hence, for $\phi \geq 0$,

$$(2.3) \qquad\qquad \langle {}^\varepsilon Y_t, \phi \rangle \leq \langle \tilde{Y}_t, \phi \rangle.$$

Similarly to Lemma 2.1, it is easy to show that

$$(2.4) \quad \langle \tilde{Y}_t, \phi \rangle = \langle f, \phi \rangle + \int_0^t \langle \tilde{Y}_r, (a\phi)'' - (b\phi)' \rangle \, dr - \int_0^t \langle \tilde{Y}_r, (c\phi)' \rangle \, dW_r.$$

Let $\phi_t$ be given by

$$(2.5) \qquad \langle f, \phi_t \rangle = \langle f, \phi \rangle + \int_0^t \langle af'' + bf', \phi_r \rangle \, dr + \int_0^t \langle cf', \phi_r \rangle \, d\tilde{W}_r,$$

where $\tilde{W}$ is an independent copy of $W$. The existence of a solution to (2.5) follows from [15]. By Itô's formula, we see that

$$e^{-\alpha \langle \tilde{Y}_t, \phi \rangle} - \int_0^t e^{-\alpha \langle \tilde{Y}_s, \phi \rangle} \left( \alpha \langle a\tilde{Y}_s'' + b\tilde{Y}_s', \phi \rangle + \frac{\alpha^2}{2} \langle c\tilde{Y}_s', \phi \rangle^2 \right) ds$$

and

$$e^{-\alpha \langle f, \phi_t \rangle} - \int_0^t e^{-\alpha \langle f, \phi_s \rangle} \left( \alpha \langle af'' + bf', \phi_s \rangle + \frac{\alpha^2}{2} \langle cf', \phi_s \rangle^2 \right) ds$$

are martingales. By a duality argument (cf. [6], page 188), we have

$$\mathbb{E} e^{-\alpha \langle \tilde{Y}_t, \phi \rangle} = \mathbb{E} e^{-\alpha \langle f, \phi_t \rangle}.$$

This implies that $\langle \tilde{Y}_t, \phi \rangle$ and $\langle f, \phi_t \rangle$ have the same distribution. Taking $f \equiv 1$ in (2.5), it is clear that

$$\int \phi_t(x) \, dx = \int \phi(x) \, dx \qquad \text{a.s.}$$

Then

$$\langle f, \phi_t \rangle \leq \|f\|_\infty \int \phi(x) \, dx \qquad \text{a.s.}$$

and hence

$$\langle \tilde{Y}_t, \phi \rangle \leq \|f\|_\infty \int \phi(x) \, dx \qquad \text{a.s.}$$

This implies the conclusion of the lemma.  □

From the proof of Lemma 2.2 , we have the following:

COROLLARY 2.3.

$$\sup_{0 \leq t \leq T} |\hat{\varepsilon y}_t - 1| \to 0 \qquad a.s.$$

as $\varepsilon \to 0$.

PROOF. From (2.4) and Condition (BC), it is easy to see that

$$\sup_{0 \leq t \leq T} \langle \tilde{Y}_t, 1 \rangle < \infty \qquad \text{a.s.}$$

The conclusion then follows from (2.3).  □

2.3. *Estimates on Sobolev norm.*  Now we give an estimate for the Sobolev norm of $^\varepsilon y_t$.

LEMMA 2.4.

$$(2.6) \qquad \mathbb{E} \sup_{0 \leq t \leq T} \|^\varepsilon y_t\|_1^4 \leq K_{2.6}.$$

PROOF. We freeze the nonlinear term and consider $^\varepsilon y_t(x)$ as the unique solution to the following linear equation:

$$z_t^\varepsilon(x) = f(x) + \int_0^t (b(x)\partial_x z_r^\varepsilon(x) + a(x)\partial_x^2 z_r^\varepsilon(x) - (T_\varepsilon{}^\varepsilon y_r^\varepsilon(x))z_r^\varepsilon(x)) \, dr$$

$$(2.7)$$

$$+ \int_0^t c(x)\partial_x z_r^\varepsilon(x) \, dW_r.$$

By [17], the solution has derivatives and their estimates depend on the bounds of $b, a, T_\varepsilon{}^\varepsilon y_r^\varepsilon, c$ and their derivatives. Since the bound of the derivative of $T_\varepsilon{}^\varepsilon y_r^\varepsilon$ may depend on $\varepsilon$, we *cannot* apply Rozovskii's estimate directly. Instead, we derive our estimate here. Note that

$$\langle z_t^\varepsilon, \phi \rangle = \langle f, \phi \rangle + \int_0^t \langle b\partial_x z_r^\varepsilon + a\partial_x^2 z_r^\varepsilon - (T_\varepsilon{}^\varepsilon y_r^\varepsilon)z_r^\varepsilon, \phi \rangle \, dr$$

$$+ \int_0^t \langle c\partial_x z_r^\varepsilon, \phi \rangle \, dW_r.$$

By Itô's formula, we have

$$\langle z_t^\varepsilon, \phi \rangle^2 = \langle f, \phi \rangle^2 + \int_0^t 2\langle z_r^\varepsilon, \phi \rangle \langle b\partial_x z_r^\varepsilon + a\partial_x^2 z_r^\varepsilon - (T_\varepsilon{}^\varepsilon y_r^\varepsilon)z_r^\varepsilon, \phi \rangle \, dr$$

$$+ \int_0^t 2\langle z_r^\varepsilon, \phi \rangle \langle c\partial_x z_r^\varepsilon, \phi \rangle \, dW_r + \int_0^t \langle c\partial_x z_r^\varepsilon, \phi \rangle^2 \, dr.$$

Adding over $\phi$ in a *complete orthonormal system* (CONS) of $H_0$, we have

$$\|z_t^\varepsilon\|_0^2 = \|f\|_0^2 + \int_0^t 2\langle z_r^\varepsilon, b\partial_x z_r^\varepsilon + a\partial_x^2 z_r^\varepsilon - (T_\varepsilon{}^\varepsilon y_r^\varepsilon)z_r^\varepsilon \rangle \, dr$$

$$+ \int_0^t 2\langle z_r^\varepsilon, c\partial_x z_r^\varepsilon \rangle \, dW_r + \int_0^t \|c\partial_x z_r^\varepsilon\|_0^2 \, dr.$$

Applying Itô's formula, we have

$$
\begin{aligned}
\|z_t^\varepsilon\|_0^4 = \|f\|_0^4 &+ \int_0^t 4\|z_r^\varepsilon\|_0^2 \langle z_r^\varepsilon, b\partial_x z_r^\varepsilon + a\partial_x^2 z_r^\varepsilon - (T_\varepsilon{}^\varepsilon y_r^\varepsilon)z_r^\varepsilon\rangle\, dr \\
&+ \int_0^t 4\|z_r^\varepsilon\|_0^2 \langle z_r^\varepsilon, c\partial_x z_r^\varepsilon\rangle\, dW_r + \int_0^t 2\|z_r^\varepsilon\|_0^2\|c\partial_x z_r^\varepsilon\|_0^2\, dr \\
&+ \int_0^t 4\langle z_r^\varepsilon, c\partial_x z_r^\varepsilon\rangle^2\, dr.
\end{aligned}
\tag{2.8}
$$

Note that the only coefficient in (2.8) which depends on $\varepsilon$ is $-(T_\varepsilon{}^\varepsilon y_r^\varepsilon)$. Since this term is negative, it can be discarded. The other terms in (2.8) can be estimated as follows: By (3.4) in [15], we have

$$
|\langle z_r^\varepsilon, b\partial_x z_r^\varepsilon\rangle| \le K_1\|z_r^\varepsilon\|_0^2 \quad \text{and} \quad |\langle z_r^\varepsilon, c\partial_x z_r^\varepsilon\rangle| \le K_2\|z_r^\varepsilon\|_0^2.
\tag{2.9}
$$

By (3.8) in [15] (with $\delta = 0$ there), we have

$$
2\langle z_r^\varepsilon, a\partial_x^2 z_r^\varepsilon\rangle + \|c\partial_x z_r^\varepsilon\|_0^2 \le K_3\|z_r^\varepsilon\|_0^2.
$$

Therefore,

$$
\|z_t^\varepsilon\|_0^4 \le \|f\|_0^4 + K_4\int_0^t \|z_r^\varepsilon\|_0^4\, dr + \int_0^t 4\|z_r^\varepsilon\|_0^2 \langle z_r^\varepsilon, c\partial_x z_r^\varepsilon\rangle\, dW_r.
$$

By the Burkholder–Davis–Gundy inequality and (2.9), we then have

$$
\begin{aligned}
\mathbb{E}\sup_{s\le t}\|z_s^\varepsilon\|_0^4 &\le \|f\|_0^4 + K_4\int_0^t \|z_r^\varepsilon\|_0^4\, dr + K_5\mathbb{E}\left(\int_0^t \|z_r^\varepsilon\|_0^4 \langle z_r^\varepsilon, c\partial_x z_r^\varepsilon\rangle^2\, dr\right)^{1/2} \\
&\le \|f\|_0^4 + K_4\int_0^t \|z_r^\varepsilon\|_0^4\, dr + K_6\mathbb{E}\left(\sup_{s\le t}\|z_s^\varepsilon\|_0^2 \left(\int_0^t \|z_r^\varepsilon\|_0^4\, dr\right)^{1/2}\right) \\
&\le \|f\|_0^4 + K_7\int_0^t \|z_r^\varepsilon\|_0^4\, dr + \tfrac{1}{2}\mathbb{E}\sup_{s\le t}\|z_s^\varepsilon\|_0^4.
\end{aligned}
$$

Therefore

$$
\mathbb{E}\sup_{s\le t}\|z_s^\varepsilon\|_0^4 \le 2\|f\|_0^4 + K_{2.10}\int_0^t \mathbb{E}\|z_r^\varepsilon\|_0^4\, dr,
\tag{2.10}
$$

where $K_{2.10}$ is a constant. Gronwall's inequality implies that

$$
\mathbb{E}\sup_{0\le t\le T}\|z_t^\varepsilon\|_0^4 \le K_{2.11}.
\tag{2.11}
$$

Let $u_r^\varepsilon = \partial_x z_r^\varepsilon$. Note that

$$
{}^\varepsilon y_r(x)\partial_x(T_\varepsilon(\hat{y}_r^\varepsilon y_r^\varepsilon)(x)) = {}^\varepsilon y_r(x)\hat{y}_r^\varepsilon T_\varepsilon u_r^\varepsilon = {}^\varepsilon y_r^\varepsilon(x) T_\varepsilon u_r^\varepsilon.
$$

Then

$$u_t^\varepsilon(x) = f'(x) + \int_0^t (c(x)\partial_x u_r^\varepsilon(x) + c'(x)u_r^\varepsilon(x)) \, dW_r$$

$$+ \int_0^t (b_1(x)\partial_x u_r^\varepsilon(x) + a(x)\partial_x^2 u_r^\varepsilon(x) + c_1^\varepsilon(x)u_r^\varepsilon(x) - {}^\varepsilon y_r^\varepsilon(x)T_\varepsilon u_r^\varepsilon(x)) \, dr,$$

where $b_1 = b + a'$, $c_1^\varepsilon = b' - T_\varepsilon{}^\varepsilon y_r^\varepsilon$. So

$$\|u_t^\varepsilon\|_0^2 = \|f'\|_0^2 + \int_0^t \|c\partial_x u_r^\varepsilon + c'u_r^\varepsilon\|_0^2 \, dr$$

$$+ \int_0^t 2\langle u_r^\varepsilon, b_1\partial_x u_r^\varepsilon + a\partial_x^2 u_r^\varepsilon + c_1^\varepsilon u_r^\varepsilon - {}^\varepsilon y_r^\varepsilon T_\varepsilon u_r^\varepsilon \rangle \, dr$$

$$+ \int_0^t 2\langle u_r^\varepsilon, c\partial_x u_r^\varepsilon + c'u_r^\varepsilon \rangle \, dW_r.$$

Note that $c_1^\varepsilon$ is bounded by a constant which does not depend on $\varepsilon$. Similar to arguments leading to (2.11), we have

(2.12) $$\mathbb{E} \sup_{0 \le t \le T} \|u_t^\varepsilon\|_0^4 \le K_{2.12}.$$

The conclusion then follows from (2.11) and (2.12). □

2.4. *Existence and uniqueness.* In this section, we prove the first part of Theorem 1.2. Let

$$z_t(x) \equiv z_t^{\varepsilon,\eta}(x) \equiv {}^\varepsilon y_t(x) - {}^\eta y_t(x).$$

Then

$$z_t(x) = \int_0^t (b(x)\partial_x z_r(x) + a(x)\partial_x^2 z_r(x) - (T_\varepsilon{}^\varepsilon y_r^\varepsilon(x){}^\varepsilon y_r(x) - T_\eta{}^\eta y_r^\eta(x){}^\eta y_r(x))) \, dr$$

$$+ \int_0^t c(x)\partial_x z_r(x) \, dW_r.$$

Note that

$$T_\varepsilon{}^\varepsilon y_r^\varepsilon{}^\varepsilon y_r - T_\eta{}^\eta y_r^\eta{}^\eta y_r = \hat{{}^\varepsilon y_r}(T_\varepsilon{}^\varepsilon y_r)z_r + \hat{{}^\varepsilon y_r}(T_\varepsilon z_r){}^\eta y_r$$

$$+ (\hat{{}^\varepsilon y_r} - \hat{{}^\eta y_r})(T_\varepsilon{}^\eta y_r){}^\eta y_r + \hat{{}^\eta y_r}(T_\varepsilon{}^\eta y_r - T_\eta{}^\eta y_r){}^\eta y_r.$$

Similarly to (2.10), we have

$$\mathbb{E} \sup_{0 \le s \le t} \|z_s\|_0^4 \le K_{2.13} \int_0^t \mathbb{E}\|z_r\|_0^4 \, dr$$

(2.13) $$+ 3\|f\|_\infty^4 \mathbb{E} \int_0^t \left( \int |T_\varepsilon{}^\eta y_r(x) - T_\eta{}^\eta y_r(x)|^2 \, dx \right)^2 \, dr$$

$$+ K_{2.13}\mathbb{E} \int_0^t |\hat{{}^\varepsilon y_r} - \hat{{}^\eta y_r}|^4 \, dr.$$

As

$$T_\varepsilon{}^\eta y_r(x) - T_\eta{}^\eta y_r(x)$$

$$= \iint_0^1 \partial_x{}^\eta y_r(x + (\theta\sqrt{\varepsilon} + (1-\theta)\sqrt{\eta})a)(\sqrt{\varepsilon} - \sqrt{\eta})a\,d\theta p(a)\,da,$$

we have, when $\varepsilon, \eta \to 0$,

$$(2.14) \qquad \int |T_\varepsilon{}^\eta y_r(x) - T_\eta{}^\eta y_r(x)|^2\,dx \le \|\partial_x{}^\eta y_r\|_0^2 (\sqrt{\varepsilon} - \sqrt{\eta})^2 \to 0,$$

where $p(a)$ is the standard normal density. By Corollary 2.3 and the dominated convergence theorem, we have

$$(2.15) \qquad \mathbb{E}\int_0^t |\hat{}^\varepsilon y_r - \hat{}^\eta y_r|^4\,dr \to 0.$$

It follows from Gronwall's inequality, (2.13)–(2.15) that

$$\mathbb{E}\sup_{0\le t\le T}\|{}^\varepsilon y_t - {}^\eta y_t\|_0^4 \to 0 \qquad \text{as } \varepsilon, \eta \to 0.$$

Hence, there exists $y_t$ s.t. ${}^\varepsilon y_t \to y_t$ in $H_0$.

Note that

$$\langle {}^\varepsilon y_t, \phi\rangle = \langle f, \phi\rangle + \int_0^t \langle {}^\varepsilon y_r, -(b\phi)' + (a\phi)'' - (T_\varepsilon{}^\varepsilon y_r^\varepsilon)\phi\rangle\,dr$$

$$+ \int_0^t \langle {}^\varepsilon y_r, -(c\phi)'\rangle\,dW_r.$$

We consider the limit of the nonlinear term only, since the other terms clearly converge to the counterpart with ${}^\varepsilon y$ replaced by $y$:

$$\mathbb{E}\left| \int_0^t \int {}^\varepsilon y_r(x)(T_\varepsilon{}^\varepsilon y_r^\varepsilon)(x)\phi(x)\,dx\,dr - \int_0^t \int y_r(x)^2\phi(x)\,dx\,dr \right|$$

$$\le \mathbb{E}\int_0^t \int |T_\varepsilon({}^\varepsilon y_r^\varepsilon - y_r)|(x){}^\varepsilon y_r(x)|\phi(x)|\,dx\,dr$$

$$+ \mathbb{E}\int_0^t \int |T_\varepsilon y_r - y_r|(x){}^\varepsilon y_r(x)|\phi(x)|\,dx\,dr$$

$$+ \mathbb{E}\int_0^t \int |{}^\varepsilon y_r - y_r|(x)y_r(x)|\phi(x)|\,dx\,dr$$

$$\to 0.$$

It is then easy to show that $y_t$ solves (1.2).

To prove the uniqueness, we assume that $y_t$ and $\tilde{y}_t$ are two solution to (1.3). Similar to (2.13), we have

$$(2.16) \qquad \mathbb{E}\sup_{s\le t}\|y_t - \tilde{y}_t\|_0^4 \le K_{2.16}\int_0^t \mathbb{E}\|y_r - \tilde{y}_r\|_0^4\,dr.$$

The uniqueness then follows from Gronwall's inequality.

LEMMA 2.5.

$$\mathbb{E} \sup_{0 \leq t \leq T} \|\partial_x y_t\|_0^4 \leq K_{2.12}.$$

PROOF. Note that

$$
\begin{aligned}
\mathbb{E} \sup_{0 \leq t \leq T} \|\partial_x y_t\|_0^4 &= \mathbb{E} \left( \sup_{0 \leq t \leq T} \sum_i \langle \partial_x y_t, \phi_i \rangle^2 \right)^2 \\
&= \mathbb{E} \left( \sup_{0 \leq t \leq T} \sum_i \langle y_t, \phi_i' \rangle^2 \right)^2 \\
&= \mathbb{E} \left( \sup_{0 \leq t \leq T} \sum_i \lim_{\varepsilon \to 0} \langle {}^{\varepsilon} y_t, \phi_i' \rangle^2 \right)^2 \\
&\leq \liminf_{\varepsilon \to 0} \mathbb{E} \left( \sup_{0 \leq t \leq T} \sum_i \langle {}^{\varepsilon} y_t, \phi_i' \rangle^2 \right)^2 \\
&= \liminf_{\varepsilon \to 0} \mathbb{E} \sup_{0 \leq t \leq T} \|\partial_x {}^{\varepsilon} y_t\|_0^4 \\
&\leq K_{2.12},
\end{aligned}
$$

where $\{\phi_i\}$ is a CONS of $H_0$.  □

2.5. *Particle representation.* In this section, we verify Theorem 1.2(ii). Let $y_t$ be the solution to (1.3) and let $Y_t(dx) = y_t(x)\,dx$. Let $(\xi_t^i, m_t^i)$ be given by (1.4) and (1.5). Denote the process given by the right-hand side of (1.6) by $\tilde{Y}_t$. Now we only need to verify that $\tilde{Y}_t$ coincides with $Y_t$. Applying Itô's formula to $m_t^i \phi(\xi_t^i)$, it is easy to show that

$$
\begin{aligned}
(2.17) \quad \langle \tilde{Y}_t, \phi \rangle &= \langle f, \phi \rangle + \int_0^t \langle \tilde{Y}_r, (a\phi)'' - (b\phi)' - y_r \phi \rangle \, dr \\
&\quad + \int_0^t \langle \tilde{Y}_r, -(c\phi)' \rangle \, dW_r.
\end{aligned}
$$

By (1.3), we see that (2.17) holds with $\tilde{Y}_t$ replaced by $Y_t$. Similar to last section, we have uniqueness for the solution of (2.17). This proves $Y_t = \tilde{Y}_t$ and hence, $Y_t$ has the particle representation given in Theorem 1.2.

**3. Wong–Zakai approximation.** In this section, we prove Theorem 1.3.

3.1. *Some estimates on $y_t^{\bar{\varepsilon}}$.* For the convenience of the reader, we state a definition and a theorem which are simplified versions of a definition on page 141 and Theorem 4.6 on page 142 in [8]. Let

$$Lu = \tilde{a}\partial_x^2 u + \tilde{b}\partial_x u + \tilde{c}u.$$

DEFINITION 3.1. A *fundamental solution* of the parabolic operator $L - \partial/\partial t$ in $\mathbb{R} \times [0, T]$ is a function $\Gamma(x, t; \xi, \tau)$ defined for all $(x, t)$ and $(\xi, \tau)$ in $\mathbb{R} \times [0, T]$, $t > \tau$, satisfying the following condition: For any continuous function $\phi(x)$ with compact support, the function

$$u(x, t) = \int_{\mathbb{R}} \Gamma(x, t; \xi, \tau)\phi(\xi)\, d\xi$$

satisfies

$$Lu - \frac{\partial u}{\partial t} = 0 \qquad \text{if } x \in \mathbb{R}, \tau < t \leq T,$$

$$u(x, t) \to \phi(x) \qquad \text{if } t \to \tau + .$$

To state the next theorem, we need the following conditions:

(A1) There is a positive constant $K$ such that

$$\tilde{a}(x, t) \geq K \qquad \text{for all } x \in \mathbb{R} \text{ and } t \in [0, T].$$

(A2) The coefficients of $L$ are bounded continuous functions in $\mathbb{R} \times [0, T]$.

(A3) The coefficients of $L$ are Hölder continuous in $x$, uniformly with respect to $(x, t)$ in compact subsets of $\mathbb{R} \times [0, T]$.

THEOREM 3.2. *Let* (A1)–(A3) *hold. Let* $g(x, t)$ *be a bounded continuous function in* $\mathbb{R} \times [0, T]$, *Hölder continuous in $x$ uniformly with respect to* $(x, t)$ *in compact subsets, and let* $\phi(x)$ *be a bounded continuous function in* $\mathbb{R}$. *Then there exists a solution of the Cauchy problem*

$$(3.1) \qquad Mu \equiv Lu(x, t) - \frac{\partial u(x, t)}{\partial t} = g(x, t) \qquad \text{in } \mathbb{R} \times [0, T]$$

*with the initial condition*

$$(3.2) \qquad u(x, 0) = \phi(x) \qquad \text{on } \mathbb{R}.$$

*The solution is given by*

$$u(x, t) = \int_{\mathbb{R}^n} \Gamma(x, t; \xi, 0)\phi(\xi)\, d\xi - \int_0^t \int_{\mathbb{R}^n} \Gamma(x, t; \xi, \tau)g(\xi, \tau)\, d\xi\, d\tau.$$

Now we come back to our equation (1.7). We shall take

$$L = \bar{a}\partial_x^2 + (\bar{b} + c\dot{W}^\varepsilon)\partial_x.$$

LEMMA 3.3.

$$(3.3) \qquad\qquad \mathbb{E}\|y_t^\varepsilon\|_0^4 \leq K_{3.3}.$$

PROOF. Given $W$, let $q^W(y,t;x,s)$ be the fundamental solution of the parabolic operator $L - \partial_t$. Then, by Theorem 3.2 and (1.7),

$$y_t^\varepsilon(x) = \int q^W(x,t;y,0) f(y)\, dy - \int_0^t \int q^W(x,t;u,s) y_s^\varepsilon(u)^2\, du\, ds$$

$$\leq \int q^W(x,t;y,0) f(y)\, dy,$$

so

$$\|y_t^\varepsilon\|_0^4 \leq \left( \int \left( \int q^W(x,t;y,0)\, dx \right) f(y)^2\, dy \right)^2$$

$$= \iiint \left( \int q^W(x_1,t;y_1,0)\, dx_1 \int q^W(x_2,t;y_2,0)\, dx_2 \right)$$

$$\times f(y_1)^2 f(y_2)^2\, dy_1\, dy_2.$$

Note that $q^W(x,t;y,0) = q^{*W}(y,0;x,t)$, $q^{*W}$ is the fundamental solution of $L^* + \partial_t$ where

$$L^*\phi = -((\bar{b} + c\dot{W}^\varepsilon)\phi)' + (\bar{a}\phi)''$$

$$= -(\bar{b}' - \bar{a}'' + c'\dot{W}^\varepsilon)\phi - (\bar{b} - 2\bar{a} + c\dot{W}^\varepsilon)\phi' + \bar{a}\phi''.$$

Let

$$d\xi_t^\varepsilon = e(\xi_t^\varepsilon)\, dB_t - (\bar{b} - 2\bar{a} + c\dot{W}_t^\varepsilon)\, dt.$$

By the Feymann–Kac formula,

$$\int q^{*W}(y,0;x,t)\, dx = E_{y,0}^W \exp\left( -\int_0^t (\bar{b}' - \bar{a}'' + c'\dot{W}_r^\varepsilon)(\xi_r^\varepsilon)\, dr \right)$$

$$\leq e^{2Kt} E_{y,0}^W \exp\left( -\int_0^t c'(\xi_r^\varepsilon)\dot{W}_r^\varepsilon\, dr \right),$$

where $E_{y,0}^W$ denotes the conditional distribution of $\xi_t^\varepsilon$ given $W$ and $\xi_0^\varepsilon = y$. Hence [assume $t = (k+1)\varepsilon$],

$$e^{-4Kt} \mathbb{E}\left( \int q^{*W}(y,0;x,t)\, dx \right)^2$$

$$\leq \mathbb{E}\left( \exp\left( -2\sum_{i=0}^{k} c'(\xi_{i\varepsilon}^\varepsilon)(W_{(i+1)\varepsilon} - W_{i\varepsilon}) \right) \right.$$

$$\times \exp\left(-2\sum_{i=0}^{k}\int_{i\varepsilon}^{(i+1)\varepsilon}\int_{i\varepsilon}^{r}((2\bar{a}-\bar{b})c'' + \bar{a}c''')(\xi_s^\varepsilon)\,ds\,\dot{W}_r^\varepsilon\,dr\right)$$

$$\times \exp\left(-2\sum_{i=0}^{k}\int_{i\varepsilon}^{(i+1)\varepsilon}\int_{i\varepsilon}^{r}c''(\xi_s^\varepsilon)(e(\xi_s^\varepsilon)\,dB_s - c(\xi_s^\varepsilon)\dot{W}_s^\varepsilon\,ds)\dot{W}_r^\varepsilon\,dr\right)\right)$$

$$\leq (I_1 I_2 I_3 I_4)^{1/4},$$

where

$$I_1 = \mathbb{E}\exp\left(-8\sum_{i=0}^{k}c'(\xi_{i\varepsilon}^\varepsilon)(W_{(i+1)\varepsilon} - W_{i\varepsilon})\right),$$

$$I_2 = \mathbb{E}\exp\left(-8\sum_{i=0}^{k}\int_{i\varepsilon}^{(i+1)\varepsilon}\int_{i\varepsilon}^{r}((2\bar{a}-\bar{b})c'' + \bar{a}c''')(\xi_s^\varepsilon)\,ds\,\dot{W}_r^\varepsilon\,dr\right),$$

$$I_3 = \mathbb{E}\exp\left(-8\sum_{i=0}^{k}\int_{i\varepsilon}^{(i+1)\varepsilon}\int_{i\varepsilon}^{r}c''(\xi_s^\varepsilon)e(\xi_s^\varepsilon)\,dB_s\,\dot{W}_r^\varepsilon\,dr\right)$$

and

$$I_4 = \mathbb{E}\exp\left(8\sum_{i=0}^{k}\int_{i\varepsilon}^{(i+1)\varepsilon}\int_{i\varepsilon}^{r}c''(\xi_s^\varepsilon)c(\xi_s^\varepsilon)\dot{W}_s^\varepsilon\,ds\,\dot{W}_r^\varepsilon\,dr\right).$$

Define $c_\varepsilon(s) = -8c'(\xi_{i\varepsilon}^\varepsilon)$ for $i\varepsilon \leq s < (i+1)\varepsilon$. Let $\tilde{P}$ be the probability measure given by

$$\frac{d\tilde{P}}{dP} = \exp\left(\int_0^t c_\varepsilon(s)\,dW_s - \frac{1}{2}\int_0^t |c_\varepsilon(s)|^2\,ds\right).$$

Then, by the Girsanov formula,

$$I_1 = \tilde{E}\exp\left(\frac{1}{2}\int_0^t |c_\varepsilon(s)|^2\,ds\right) \leq \exp(32\|c'\|_\infty^2 t),$$

where $\tilde{E}$ denotes the expectation under the measure $\tilde{P}$. Note that for $\varepsilon$ small, more precisely, for

$$\varepsilon < \min((4\|(2\bar{a}-\bar{b})c'' + \bar{a}c'''\|_\infty)^{-1/2}, (8\|ec'\|_\infty)^{-1}),$$

we have

$$I_2 \leq \mathbb{E}\exp\left(4\|(2\bar{a}-\bar{b})c'' + \bar{a}c'''\|_\infty \sum_{i=0}^{k}\varepsilon|W_{(i+1)\varepsilon} - W_{i\varepsilon}|\right)$$

$$\leq \mathbb{E}\exp\left(2\|(2\bar{a}-\bar{b})c'' + \bar{a}c'''\|_\infty\left(t + \varepsilon\sum_{i=0}^{k}|W_{(i+1)\varepsilon} - W_{i\varepsilon}|^2\right)\right)$$

$$\leq \exp(2\|(2\bar{a}-\bar{b})c'' + \bar{a}c'''\|_\infty t)(1 - 4\|(2\bar{a}-\bar{b})c'' + \bar{a}c'''\|_\infty \varepsilon^2)^{-k/2}$$

$$\leq \exp(10\|(2\bar{a}-\bar{b})c'' + \bar{a}c'''\|_\infty t),$$

$$I_3 = \mathbb{E}\mathbb{E}^W \exp\left(-8 \sum_{i=0}^{k} \int_{i\varepsilon}^{(i+1)\varepsilon} c''(\xi_s^\varepsilon) e(\xi_s^\varepsilon)\varepsilon^{-1} \right.$$
$$\left. \times ((i+1)\varepsilon - s)(W_{(i+1)\varepsilon} - W_{i\varepsilon})\,dB_s\right)$$

$$\leq \mathbb{E}\exp\left(32 \sum_{i=0}^{k} \int_{i\varepsilon}^{(i+1)\varepsilon} c''(\xi_s^\varepsilon)^2 e(\xi_s^\varepsilon)^2 (W_{(i+1)\varepsilon} - W_{i\varepsilon})^2\,ds\right)$$

$$\leq \mathbb{E}\exp\left(32\|ec''\|_\infty^2 \sum_{i=0}^{k} (W_{(i+1)\varepsilon} - W_{i\varepsilon})^2 \varepsilon\right)$$

$$\leq \prod_{i=0}^{k} (1 - 64\|ec''\|_\infty^2 \varepsilon^2)^{-1/2}$$

$$\leq \exp(32\|ec''\|_\infty^2 \varepsilon t)$$

and

$$I_4 \leq \mathbb{E}\exp\left(8 \sum_{i=0}^{k} \|cc''\|_\infty (W_{(i+1)\varepsilon} - W_{i\varepsilon})^2\right)$$

$$\leq \exp(32\|cc''\|_\infty t).$$

The conclusion then follows easily. $\square$

We now turn to the estimation on the norm of $\partial_x y_t^\varepsilon$.

LEMMA 3.4. *Suppose that* $\{N(x) : x \in \mathbb{R}\}$ *is a random field such that* $\exists\, \alpha > 0$, $p > 1$,

$$\mathbb{E}(|N(x) - N(y)|^p) \leq K|x - y|^{1+\alpha}.$$

*Then for any* $\lambda > 0$,

$$\mathbb{E}\sup_{x \in \mathbb{R}}(|N(x)|^p e^{-\lambda|x|}) < \infty.$$

PROOF. It follows from Theorem 4 in [10] that, for any $I_n = [n, n+1]$,

$$\left(\mathbb{E}\sup_{x,y \in I_n} |N(x) - N(y)|^p\right)^{1/p} \leq C \int_0^1 \frac{Ku^{(1+\alpha)/p}}{u^{1+1/p}}\,du$$

$$\leq \frac{CKp}{\alpha} \equiv K_1.$$

Note that

$$|N(x) - N(0)|^p e^{-\lambda|x|}$$

$$\leq \left( \sum_n \sup_{y,z \in I_n} |N(y) - N(z)| e^{-\lambda|n|/p} \right)^p$$

$$\leq (2(1 - e^{-\lambda/p}))^{(1-p)/p} \sum_n \sup_{y,z \in I_n} |N(y) - N(z)|^p e^{-\lambda|n|/p}.$$

Hence

$$\mathbb{E} \sup_{x \in \mathbb{R}} (|N(x) - N(0)|^p e^{-\lambda|x|})$$

$$\leq (2(1 - e^{-\lambda/p}))^{(1-p)/p} \sum_n \mathbb{E} \sup_{y,z \in I_n} |N(y) - N(z)|^p e^{-\lambda|n|/p}$$

$$\leq (2(1 - e^{-\lambda/p}))^{(1-p)/p} \sum_n K_1^p e^{-\lambda|n|/p}$$

$$\leq K_1^p (2(1 - e^{-\lambda/p}))^{(1-2p)/p} < \infty.$$

The conclusion of the lemma then follows easily.  $\square$

LEMMA 3.5.

$$(3.4) \qquad\qquad \mathbb{E}\|\partial_x y_t^\varepsilon\|_0^4 \leq K_{3.4}.$$

PROOF.   Note that

$$\partial_x y_t^\varepsilon = f' + \int_0^t ((\bar{b}' - 2y_r^\varepsilon + c'\dot{W}_r^\varepsilon)\partial_x y_r^\varepsilon$$

$$+ (\bar{b} + \bar{a}' + c\dot{W}_r^\varepsilon)\partial_x^2 y_r^\varepsilon + \bar{a}\partial_x^3 y_r^\varepsilon) \, dr.$$

Let $q_1^W$ be the fundamental solution of $L_1 - \partial_t$, where

$$L_1\phi = \bar{a}\phi'' + (\bar{b} + \bar{a}' + c\dot{W}_r^\varepsilon)\phi' + (\bar{b}' - 2y_r^\varepsilon + c'\dot{W}_r^\varepsilon)\phi.$$

Then

$$\partial_x y_t^\varepsilon = \int q_1^W(x, t; y, 0) f'(y) \, dy.$$

Note that

$$L_1^*\phi = (\bar{a}\phi)'' - ((\bar{b} + \bar{a}' + c\dot{W}_r^\varepsilon)\phi)' + (\bar{b}' - 2y_r^\varepsilon + c'\dot{W}_r^\varepsilon)\phi$$

$$= \bar{a}\phi'' + (\bar{a}' - \bar{b} - c\dot{W}_r^\varepsilon)\phi' - 2y_r^\varepsilon\phi.$$

Similarly to Lemma 3.3, we have, for any $\lambda$ and $p > 1$,

$$(3.5) \qquad\qquad \mathbb{E}\left( \int e^{\lambda|x|} q_1^W(x, t; y, 0) \, dx \right)^p \leq K_1$$

and

$$\int q_1^W(x,t;y,0)\,dy = \mathbb{E}_{0,x}^W \exp\left(\int_0^t (\bar{b}' - 2y_r^\varepsilon + c'\dot{W}_r^\varepsilon)(\eta_r^{\varepsilon,x})\,dr\right)$$

$$\leq e^{\|\bar{b}'\|_\infty} \mathbb{E}_{0,x}^W \exp\left(\sum_{i=0}^k \int_{i\varepsilon}^{(i+1)\varepsilon} c'(\eta_r^{\varepsilon,x})\,dr\,\dot{W}_{i\varepsilon}^\varepsilon\right),$$

where $\eta_t^{\varepsilon,x}$, with initial $x$, solves

$$d\eta_t^{\varepsilon,x} = (\bar{b}+\bar{a})(\eta_t^{\varepsilon,x})\,dt + c(\eta_t^{\varepsilon,x})\dot{W}_t^\varepsilon\,dt + e(\eta_t^{\varepsilon,x})\,dB_t.$$

Note that for $i\varepsilon \leq r \leq (i+1)\varepsilon$,

$$c'(\eta_r^{\varepsilon,x}) = c'(\eta_{i\varepsilon}^{\varepsilon,x}) + \int_{i\varepsilon}^r c''(\eta_s^{\varepsilon,x})e(\eta_s^{\varepsilon,x})\,dB_s$$

$$+ \int_{i\varepsilon}^r \left((\bar{a}+\bar{b})c'' + \frac{e^2}{2}c'''\right)ds + \int_{i\varepsilon}^r cc''\,ds\,\dot{W}_{i\varepsilon}^\varepsilon.$$

As

$$\left|\sum_{i=0}^k \int_{i\varepsilon}^{(i+1)\varepsilon} \left(\int_{i\varepsilon}^r \left((\bar{a}-\bar{b})c'' + \frac{e^2}{2}c'''\right)ds + \int_{i\varepsilon}^r cc''\,ds\,\dot{W}_{i\varepsilon}^\varepsilon\right)dr\,\dot{W}_{i\varepsilon}^\varepsilon\right|$$

$$\leq \sum_{i=0}^k \left\|(\bar{a}-\bar{b})c'' + \frac{e^2}{2}c'''\right\|_\infty \varepsilon(W_{(i+1)\varepsilon} - W_{i\varepsilon})$$

$$+ \sum_{i=0}^k \|cc''\|_\infty (W_{(i+1)\varepsilon} - W_{i\varepsilon})^2$$

$$\equiv \frac{1}{4}\log M(W) - \|\bar{b}'\|_\infty,$$

we have

$$\left(\int q_1^W(x,t;y,0)\,dy\right)^4 \leq M(W)N_1(x,W)N_2(x,W),$$

where

$$N_1(x,W) = \mathbb{E}_{0,x}^W \exp\left(-4\sum_{i=0}^k c'(\eta_{i\varepsilon}^{\varepsilon,x})(W_{(i+1)\varepsilon} - W_{i\varepsilon})\right)$$

and

$$N_2(x,W) = \mathbb{E}_{0,x}^W \exp\left(-4\sum_{i=0}^k \int_{i\varepsilon}^{(i+1)\varepsilon} \int_{i\varepsilon}^r c''(\eta_s^{\varepsilon,x})e(\eta_s^{\varepsilon,x})\,dB_s\,\dot{W}_r^\varepsilon\,dr\right).$$

By arguments similar to Lemma 3.3, it is easy to see that $M(W)$ has finite moments.

First take $\mathbb{E}^W$ and then take expectation with respect to $W$, for $t \in [i\varepsilon, (i+1)\varepsilon]$ and even integer $p$; we have

$$\mathbb{E}|\eta_t^{\varepsilon,x} - \eta_t^{\varepsilon,y}|^p$$

$$\leq \mathbb{E}|\eta_{i\varepsilon}^{\varepsilon,x} - \eta_{i\varepsilon}^{\varepsilon,y}|^p + \int_{i\varepsilon}^t K_1 \mathbb{E}|\eta_s^{\varepsilon,x} - \eta_s^{\varepsilon,y}|^p \, ds$$

$$+ p\mathbb{E}\int_{i\varepsilon}^t (c(\eta_s^{\varepsilon,x}) - c(\eta_s^{\varepsilon,y}))(\eta_s^{\varepsilon,x} - \eta_s^{\varepsilon,y})^{p-1} \dot{W}_s^\varepsilon \, ds$$

$$\leq \mathbb{E}|\eta_{i\varepsilon}^{\varepsilon,x} - \eta_{i\varepsilon}^{\varepsilon,y}|^p + \int_{i\varepsilon}^t K_2 \mathbb{E}|\eta_s^{\varepsilon,x} - \eta_s^{\varepsilon,y}|^p \, ds$$

$$+ p(p-1)\mathbb{E}\int_{i\varepsilon}^t \int_{i\varepsilon}^s (c(\eta_r^{\varepsilon,x}) - c(\eta_r^{\varepsilon,y}))^2$$

$$\times (\eta_r^{\varepsilon,x} - \eta_r^{\varepsilon,y})^{p-2} \, dr \, ds \varepsilon^{-2}(W_{(i+1)\varepsilon} - W_{i\varepsilon})^2$$

$$\leq (1 + K_3\varepsilon)\mathbb{E}|\eta_{i\varepsilon}^{\varepsilon,x} - \eta_{i\varepsilon}^{\varepsilon,y}|^p.$$

By induction we have

$$\mathbb{E}|\eta_t^{\varepsilon,x} - \eta_t^{\varepsilon,y}|^p \leq K_2 |x-y|^p.$$

Therefore

$$\mathbb{E}|N_i(x,W) - N_i(y,W)|^p \leq K_3 |x-y|^{p/2}, \qquad i=1,2.$$

By Lemma 3.4 we have

$$\mathbb{E}\sup_x |N_i(x,W)|^p e^{-\lambda|x|} \leq K_4.$$

Therefore

$$\mathbb{E}\sup_x \left( \int q_1^W(x,t;y,0) \, dy e^{-\lambda|x|} \right)^4$$

$$\leq \mathbb{E}\left( M(W) \sup_x N_1(x,W) e^{-2\lambda|x|} \sup_x N_2(x,W) e^{-2\lambda|x|} \right)$$

$$\leq K_5.$$

Note that

$$\int (\partial_x y_t^\varepsilon)(x)^2 \, dx$$

$$\leq \int \left( \int q_1^W(x,t;y,0)|f'(y)| \, dy \int q_1^W(x,t;y,0) \, dy \right) dx \|f'\|_\infty$$

$$\leq \int \left( \int e^{\lambda|x|} q_1^W(x,t;y,0) \, dx \right) |f'(y)| \, dy$$

$$\times \sup_x \int q_1^W(x,t;y,0) \, dy e^{-\lambda|x|} \|f'\|_\infty.$$

Hence

$$(\mathbb{E}\|\partial_x y_t^\varepsilon\|_0^4)^2 \leq \|f'\|_\infty^4 \mathbb{E}\left(\int \left(\int e^{\lambda|x|} q_1^W(x,t;y,0)\,dx\right)|f'(y)|\,dy\right)^4$$

$$\times \mathbb{E}\left(\sup_x \int q_1^W(x,t;y,0)\,dy e^{-\lambda|x|}\right)^4$$

$$\leq \|f'\|_\infty^4 \mathbb{E}\int\left(\int e^{\lambda|x|} q_1^W(x,t;y,0)\,dx\right)^4$$

$$\times |f'(y)|^2\,dy\left(\int |f'(y)|^{2/3}\,dy\right)^3 K_5$$

$$\leq \|f'\|_\infty^4 K_1 \|f'\|_0^2\left(\int |f'(y)|^{2/3}\,dy\right)^3 K_5 < \infty.$$

This proves the conclusion of the lemma. $\square$

COROLLARY 3.6.  (i) *For any $\alpha \geq 0$ and $p \geq 0$, we have*

$$(3.6) \qquad \mathbb{E}\left|\int_\mathbb{R} |\partial_x y_t^\varepsilon(x)|^{1+\alpha}\,dx\right|^p \leq K_{3.6}.$$

(ii)

$$(3.7) \qquad \mathbb{E}\|\partial_x^2 y_t^\varepsilon\|_0^4 \leq K_{3.7}.$$

PROOF.  The proof of Lemma 3.5 can be modified to verify (i). Part (ii) follows from the same proof as well; note that (i) implies $\mathbb{E}\|(\partial_x y_t^\varepsilon)^2\|_0^4 \leq K_{3.6}$. $\square$

3.2. *Proof of Theorem 1.3.*  Now we prove Theorem 1.3. In this proof, the quantity $\langle \partial_x^2 z_r^\varepsilon, f\rangle$ for $f$ smooth is understood as $\langle z_r^\varepsilon, \partial_x^2 f\rangle$.

To make use of Itô's formula, we need that $y_t^\varepsilon$ is adapted. We shall use $y_{t-\varepsilon}^\varepsilon$ to replace $y_t^\varepsilon$. However, for simplicity of notation, we still use $y_t^\varepsilon$.

Let $z_t^\varepsilon = y_t^\varepsilon - y_t$. Then

$$\langle z_t^\varepsilon, \phi\rangle = \int_0^t \langle b\partial_x z_r^\varepsilon + a\partial_x^2 z_r^\varepsilon - (y_r^\varepsilon + y_r)z_r^\varepsilon, \phi\rangle\,dr$$

$$+ \int_0^t \langle c\partial_x y_r^\varepsilon, \phi\rangle \dot{W}_{r-\varepsilon}^\varepsilon\,dr$$

$$- \int_0^t \langle c\partial_x y_r, \phi\rangle\,dW_r - \int_0^t \langle \tfrac{1}{2}cc'\partial_x y_r^\varepsilon + \tfrac{1}{2}c^2\partial_x^2 y_r^\varepsilon, \phi\rangle\,dr.$$

By Itô's formula, we have

$$\langle z_t^\varepsilon, \phi\rangle^2 = \int_0^t 2\langle z_r^\varepsilon, \phi\rangle\langle b\partial_x z_r^\varepsilon + a\partial_x^2 z_r^\varepsilon - (y_r^\varepsilon + y_r)z_r^\varepsilon, \phi\rangle\,dr$$

$$+ \int_0^t 2\langle z_r^\varepsilon, \phi\rangle\langle c\partial_x y_r^\varepsilon, \phi\rangle \dot{W}_{r-\varepsilon}^\varepsilon \, dr - \int_0^t 2\langle z_r^\varepsilon, \phi\rangle\langle c\partial_x y_r, \phi\rangle \, dW_r$$

$$- \int_0^t \langle z_r^\varepsilon, \phi\rangle\langle cc'\partial_x y_r^\varepsilon + c^2\partial_x^2 y_r^\varepsilon, \phi\rangle \, dr$$

$$+ \int_0^t \langle c\partial_x y_r, \phi\rangle^2 \, dr.$$

Adding over $\phi$ in a CONS of $H_0$, we have

$$
\begin{aligned}
(3.8) \quad \|z_t^\varepsilon\|_0^2 = {} & \int_0^t 2\langle z_r^\varepsilon, b\partial_x z_r^\varepsilon + a\partial_x^2 z_r^\varepsilon - (y_r^\varepsilon + y_r)z_r^\varepsilon\rangle \, dr \\
& + \int_0^t 2\langle z_r^\varepsilon, c\partial_x y_r^\varepsilon\rangle \dot{W}_{r-\varepsilon}^\varepsilon \, dr - \int_0^t 2\langle z_r^\varepsilon, c\partial_x y_r\rangle \, dW_r \\
& - \int_0^t \langle z_r^\varepsilon, cc'\partial_x y_r^\varepsilon + c^2\partial_x^2 y_r^\varepsilon\rangle \, dr \\
& + \int_0^t \|c\partial_x y_r\|_0^2 \, dr.
\end{aligned}
$$

We now estimate the second term on the right-hand side of (3.8). For $(i-1)\varepsilon \le r < i\varepsilon$, note that

$$
\begin{aligned}
\langle z_r^\varepsilon, \phi\rangle = {} & \langle z_{(i-1)\varepsilon}^\varepsilon, \phi\rangle + \int_{(i-1)\varepsilon}^r \langle b\partial_x z_s^\varepsilon + a\partial_x^2 z_s^\varepsilon - (y_s^\varepsilon + y_s)z_s^\varepsilon, \phi\rangle \, ds \\
& + \int_{(i-1)\varepsilon}^r \langle c\partial_x y_s^\varepsilon, \phi\rangle \dot{W}_{s-\varepsilon}^\varepsilon \, ds \\
& - \int_{(i-1)\varepsilon}^r \langle c\partial_x y_s, \phi\rangle \, dW_s - \int_{(i-1)\varepsilon}^r \langle \tfrac{1}{2}cc'\partial_x y_s^\varepsilon + \tfrac{1}{2}c^2\partial_x^2 y_s^\varepsilon, \phi\rangle \, ds
\end{aligned}
$$

and

$$
\begin{aligned}
\langle c\partial_x y_r^\varepsilon, \phi\rangle = {} & \langle c\partial_x y_{(i-1)\varepsilon}^\varepsilon, \phi\rangle + \int_{(i-1)\varepsilon}^r \langle c\partial_x(\bar{b}\partial_x y_s^\varepsilon + \bar{a}\partial_x^2 y_s^\varepsilon - (y_s^\varepsilon)^2), \phi\rangle \, ds \\
& + \int_{(i-1)\varepsilon}^r \langle c\partial_x(c\partial_x y_s^\varepsilon), \phi\rangle \dot{W}_{s-\varepsilon}^\varepsilon \, ds.
\end{aligned}
$$

Similarly to (3.8), we have

$$
\begin{aligned}
(3.9) \quad & \langle z_r^\varepsilon, c\partial_x y_r^\varepsilon\rangle - \langle z_{(i-1)\varepsilon}^\varepsilon, c\partial_x y_{(i-1)\varepsilon}^\varepsilon\rangle \\
& = \int_{(i-1)\varepsilon}^r \langle c\partial_x y_s^\varepsilon, b\partial_x z_s^\varepsilon + a\partial_x^2 z_s^\varepsilon - (y_s^\varepsilon + y_s)z_s^\varepsilon\rangle \, ds \\
& \quad + \int_{(i-1)\varepsilon}^r \|c\partial_x y_s^\varepsilon\|_0^2 \dot{W}_{s-\varepsilon}^\varepsilon \, ds - \int_{(i-1)\varepsilon}^r \langle c\partial_x y_s^\varepsilon, c\partial_x y_s\rangle \, dW_s
\end{aligned}
$$

$$- \int_{(i-1)\varepsilon}^{r} \langle c\partial_x y_s^\varepsilon, \tfrac{1}{2}cc'\partial_x y_s^\varepsilon + \tfrac{1}{2}c^2\partial_x^2 y_s^\varepsilon\rangle\, ds$$

$$+ \int_{(i-1)\varepsilon}^{r} \langle z_s^\varepsilon, c\partial_x(\bar b\partial_x y_s^\varepsilon + \bar a\partial_x^2 y_s^\varepsilon - (y_s^\varepsilon)^2)\rangle\, ds$$

$$+ \int_{(i-1)\varepsilon}^{r} \langle z_s^\varepsilon, c\partial_x(c\partial_x y_s^\varepsilon)\rangle \dot W_{s-\varepsilon}^\varepsilon\, ds.$$

Let $t = k\varepsilon$. Then

$$\mathbb{E}\int_0^t 2\langle z_r^\varepsilon, c\partial_x y_r^\varepsilon\rangle \dot W_{r-\varepsilon}^\varepsilon\, dr$$

$$(3.10) \qquad = \mathbb{E}\sum_{i=0}^{k-1} \int_{i\varepsilon}^{(i+1)\varepsilon} 2\langle z_r^\varepsilon, c\partial_x y_r^\varepsilon\rangle \dot W_{r-\varepsilon}^\varepsilon\, dr$$

$$= \mathbb{E}\sum_{i=0}^{k-1} 2\int_{i\varepsilon}^{(i+1)\varepsilon} (\langle z_r^\varepsilon, c\partial_x y_r^\varepsilon\rangle - \langle z_{(i-1)\varepsilon}^\varepsilon, c\partial_x y_{(i-1)\varepsilon}^\varepsilon\rangle)\dot W_{r-\varepsilon}^\varepsilon\, dr.$$

Apply (3.9) to (3.10). We only consider the second, third and sixth terms in (3.9) since it is easy to verify that the other terms result in quantities bounded by $K\sqrt\varepsilon$. Note that

$$\mathbb{E}\sum_{i=0}^{k-1} 2\int_{i\varepsilon}^{(i+1)\varepsilon} \int_{(i-1)\varepsilon}^{r} \|c\partial_x y_s^\varepsilon\|_0^2 \dot W_{s-\varepsilon}^\varepsilon\, ds\, \dot W_{r-\varepsilon}^\varepsilon\, dr$$

$$\approx \mathbb{E}\sum_{i=0}^{k-1} 2\int_{i\varepsilon}^{(i+1)\varepsilon} \int_{(i-1)\varepsilon}^{i\varepsilon} \|c\partial_x y_{(i-2)\varepsilon}^\varepsilon\|_0^2 \dot W_{s-\varepsilon}^\varepsilon\, ds\, \dot W_{r-\varepsilon}^\varepsilon\, dr$$

$$(3.11) \qquad + \mathbb{E}\sum_{i=0}^{k-1} 2\int_{i\varepsilon}^{(i+1)\varepsilon} \int_{i\varepsilon}^{r} \|c\partial_x y_{(i-2)\varepsilon}^\varepsilon\|_0^2\, ds\, dr\, \varepsilon^{-2}(W_{i\varepsilon} - W_{(i-1)\varepsilon})^2$$

$$= \sum_{i=0}^{k-1} \varepsilon^2 \mathbb{E}\|c\partial_x y_{(i-2)\varepsilon}^\varepsilon\|_0^2 \varepsilon^{-2}\varepsilon$$

$$\approx \mathbb{E}\int_0^t \|c\partial_x y_r^\varepsilon\|_0^2\, dr,$$

where $x \approx y$ means that $|x - y| \le K\sqrt\varepsilon$. Similarly,

$$\mathbb{E}\sum_{i=0}^{k-1} 2\int_{i\varepsilon}^{(i+1)\varepsilon} \int_{(i-1)\varepsilon}^{r} \langle z_s^\varepsilon, c\partial_x(c\partial_x y_s^\varepsilon)\rangle \dot W_{s-\varepsilon}^\varepsilon\, ds\, \dot W_{r-\varepsilon}^\varepsilon\, dr$$

$$(3.12) \qquad \approx \mathbb{E}\int_0^t \langle z_r^\varepsilon, c\partial_x(c\partial_x y_r^\varepsilon)\rangle\, dr.$$

Note that

$$\mathbb{E} \sum_{i=0}^{k-1} 2 \int_{i\varepsilon}^{(i+1)\varepsilon} \int_{(i-1)\varepsilon}^{r} \langle c\partial_x y_s^\varepsilon, c\partial_x y_s \rangle \, dW_s \dot{W}_{r-\varepsilon}^\varepsilon \, dr$$

$$(3.13) \qquad \approx \mathbb{E} \sum_{i=0}^{k-1} 2 \int_{i\varepsilon}^{(i+1)\varepsilon} \int_{(i-1)\varepsilon}^{r} \langle c\partial_x y_{(i-2)\varepsilon}^\varepsilon, c\partial_x y_{(i-2)\varepsilon} \rangle \, dW_s \dot{W}_{r-\varepsilon}^\varepsilon \, dr$$

$$\approx 2\mathbb{E} \int_0^t \langle c\partial_x y_r^\varepsilon, c\partial_x y_r \rangle \, dr.$$

By (3.8) and (3.10)–(3.13), we have

$$\mathbb{E}\|z_t^\varepsilon\|_0^2 \le K_1 \int_0^t \mathbb{E}\|z_s^\varepsilon\|_0^2 \, ds + K_2\sqrt{\varepsilon}.$$

Gronwall's inequality then implies the conclusion of the theorem.

### 4. Log-Laplace transform of $X_t$.
In this section we prove Theorem 1.4. Since $X^\varepsilon$ solves the CMP defined in Section 1.2, it is easy to show that

$$\mathbb{E} \sup_{s \le t} \langle X_s^\varepsilon, 1 \rangle^4 \le K_1.$$

For any $\phi \in C_c(\mathbb{R})$, it is then easy to show that

$$\mathbb{E} \langle X_t^\varepsilon - X_s^\varepsilon, \phi \rangle^4 \le K_2 |t-s|^2.$$

This implies the tightness of $\{X^\varepsilon\}$ in $C([0,\infty), \mathcal{M}_F(\bar{\mathbb{R}}))$. It is then easy to verify that any one of the limit points solves the MP. The uniqueness for the solution to the MP implies the weak convergence of $X^\varepsilon$.

Now we prove (1.8). First we assume that $\mu \in H_0$ and fix $X_0^\varepsilon = \mu$. Let $\psi$ be a bounded continuous function on $C([0,t], \mathbb{R})$. Then

$$\mathbb{E}(\exp(-\langle X_t, f \rangle)\psi(W)) = \lim_{\varepsilon \to 0} \mathbb{E}(\exp(-\langle X_t^\varepsilon, f \rangle)\psi(W))$$

$$= \lim_{\varepsilon \to 0} \mathbb{E}(\exp(-\langle \mu, y_{0,t}^\varepsilon \rangle)\psi(W))$$

$$= \mathbb{E}(\exp(-\langle \mu, y_{0,t} \rangle)\psi(W)).$$

For general $\mu$, we take $\mu^\varepsilon \in H_0$ converging to $\mu$ in $M_F(\mathbb{R})$. Denote the solution of the MP with $\mu$ replaced by $\mu^\varepsilon$ by $X^{(\varepsilon)}$. Then

$$\mathbb{E}(\exp(-\langle X_t, f \rangle)\psi(W)) = \lim_{\varepsilon \to 0} \mathbb{E}(\exp(-\langle X_t^{(\varepsilon)}, f \rangle)\psi(W))$$

$$= \lim_{\varepsilon \to 0} \mathbb{E}(\exp(-\langle \mu^\varepsilon, y_{0,t} \rangle)\psi(W))$$

$$= \mathbb{E}(\exp(-\langle \mu, y_{0,t} \rangle)\psi(W)),$$

where the last equation follows since $y_{0,t}$ is bounded and continuous.

**5. Moments of $X_t$.** In this section, we prove Theorem 1.5. Let $y_t^\alpha$ be the solution of

$$
\begin{aligned}
y_t^\alpha(x) = \alpha f(x) + \int_0^t (b(x)\partial_x y_r^\alpha(x) + a(x)\partial_x^2 y_r^\alpha(x) - y_r^\alpha(x)^2)\, dr \\
+ \int_0^t c(x)\partial_x y_r^\alpha(x)\, dW_r.
\end{aligned}
\tag{5.1}
$$

Let $z_t$ and $h_t$ be solutions to

$$
\begin{aligned}
z_t(x) = f(x) + \int_0^t (b(x)\partial_x z_r(x) + a(x)\partial_x^2 z_r(x))\, dr \\
+ \int_0^t c(x)\partial_x z_r(x)\, dW_r
\end{aligned}
\tag{5.2}
$$

and

$$
\begin{aligned}
h_t(x) = \int_0^t (b(x)\partial_x h_r(x) + a(x)\partial_x^2 h_r(x) - 2z_r(x)^2)\, dr \\
+ \int_0^t c(x)\partial_x h_r(x)\, dW_r.
\end{aligned}
\tag{5.3}
$$

Define $z_t^\alpha = \alpha^{-1} y_t^\alpha - z_t$. Then

$$
\begin{aligned}
z_t^\alpha(x) = \int_0^t (b(x)\partial_x z_r^\alpha(x) + a(x)\partial_x^2 z_r^\alpha(x))\, dr \\
+ \int_0^t c(x)\partial_x z_r^\alpha(x)\, dW_r - \int_0^t \alpha^{-1} y_r^\alpha(x)^2\, dr.
\end{aligned}
$$

Similarly to arguments in previous sections, we have

$$
\mathbb{E}\|z_t^\alpha\|_0^2 \to 0 \qquad \text{as } \alpha \to 0.
$$

Define $h_t^\alpha = \alpha^{-2}(y_t^{2\alpha} - 2y_t^\alpha) - h_t$. Then

$$
\begin{aligned}
h_t^\alpha(x) = \int_0^t (b(x)\partial_x h_r^\alpha(x) + a(x)\partial_x^2 h_r^\alpha(x))\, dr \\
+ \int_0^t c(x)\partial_x h_r^\alpha(x)\, dW_r \\
- \int_0^t (\alpha^{-2}(y_r^{2\alpha}(x)^2 - 2y_r^\alpha(x)^2) - 2z_r(x)^2)\, dr.
\end{aligned}
$$

Note that $|y_r^\alpha(x)| \le \alpha\|f\|_\infty$ and $|z_r(x)| \le \|f\|_\infty$. Hence

$$
\mathbb{E}\int (\alpha^{-2}(y_r^{2\alpha}(x)^2 - 2y_r^\alpha(x)^2) - 2z_r(x)^2)^2\, dx
$$

$$= \mathbb{E} \int \left( 4 \left( \frac{y_r^{2\alpha}(x)}{2\alpha} - z_r(x) \right)^2 - 2 \left( \frac{y_r^{\alpha}(x)}{\alpha} - z_r(x) \right)^2 \right.$$
$$\left. + 4 z_r(x) \frac{y_r^{2\alpha}(x) - y_r^{\alpha}(x) - \alpha z_r(x)}{\alpha} \right)^2 dx$$
$$\to 0.$$

Similarly to above, we have

$$\mathbb{E} \| h_t^{\alpha} \|_0^2 \to 0 \qquad \text{as } \alpha \to 0.$$

Therefore $z_t = \partial_{\alpha} y_t^{\alpha}|_{\alpha=0}$ and $h_t = \partial_{\alpha}^2 y_t^{\alpha}|_{\alpha=0}$.

Note that

$$\mathbb{E}(\langle X_t, f \rangle | W) = \langle \mu, z_t \rangle$$

and

$$\mathbb{E}(\langle X_t, f \rangle^2 | W) = \langle \mu, z_t \rangle^2 - \langle \mu, h_t \rangle.$$

Taking expectations on both sides of (5.2), we have

$$\mathbb{E} z_t(x) = f(x) + \mathbb{E} \int_0^t (b(x) \partial_x z_r(x) + a(x) \partial_x^2 z_r(x)) \, dr,$$

and hence, (1.9) holds.

Applying Itô's formula to (5.2), we have

$$\mathbb{E} z_t(x_1) z_t(x_2)$$
$$= f(x_1) f(x_2) + \mathbb{E} \int (b(x_1) \partial_{x_1} z_r(x_1) z_r(x_2) + b(x_2) \partial_{x_2} z(x_1) z_r(x_2)$$
$$+ a(x_1) \partial_{x_1}^2 z(x_1) z_r(x_2) + a(x_2) \partial_{x_2}^2 z(x_1) z_r(x_2)$$
$$+ c(x_1) c(x_2) \partial_{x_1} \partial_{x_2} z(x_1) z_r(x_2)) \, dr.$$

Hence

$$(5.4) \qquad \mathbb{E} z_t(x_1) z_t(x_2) = \iint f(y_1) f(y_2) q(t, (x_1, x_2), (y_1, y_2)) \, dy_1 \, dy_2.$$

Taking the expectation on both sides of (5.3), we have

$$(5.5) \qquad \mathbb{E} h_t(x) = \mathbb{E} \int_0^t (b(x) \partial_x h_r(x) + a(x) \partial_x^2 h_r(x) - 2 z_r(x)^2) \, dr.$$

Hence, making use of (5.4) and solving (5.5), we obtain

$$\mathbb{E} h_t(x) = -2 \int_0^t \int p(t - s, x, y) \mathbb{E} z_s(y)^2 \, dy \, ds$$

$$= -2 \int_0^t \int p(t-s,x,y) \iint f(z_1) f(z_2)$$

$$\times \, q(s,(y,y),(z_1,z_2)) \, dz_1 \, dz_2 \, dy \, ds.$$

This proves (1.10).

**Acknowledgments.** I am grateful to Robert Adler for his suggestion which led to Theorem 1.4. I would like to thank the referee for pointing out two gaps overlooked by me in the original manuscript. I would like to thank both the Associate Editor and the referee for helpful suggestions. I thank Klaus Fleischmann for going over the referee's report with me and for a careful reading of the manuscript.

DEPARTMENT OF MATHEMATICS
UNIVERSITY OF TENNESSEE
KNOXVILLE, TENNESSEE 37996-1300
USA
E-MAIL: jxiong@math.utk.edu